\documentclass [12pt,a4paper,reqno]{amsart}
\usepackage{amsmath,amsthm,amscd}
\usepackage{amsfonts,amssymb}
\allowdisplaybreaks

\newcommand{\be}{\begin{equation}}
\newcommand{\ef}{\end{equation}}
\chardef\bslash=`\\ % p. 424, TeXbook
\newtheorem*{thm*}{Theorem}

\theoremstyle{definition}

\newtheorem*{remark*}{Remarks}
\newtheorem*{defn*}{Definition}

\theoremstyle{remark}
\newcommand{\G}{\Gamma}
\newcommand{\wt}{\widetilde}
\newcommand{\wh}{\widehat}
\newcommand{\bk}{\bigskip}
 \renewcommand{\sectionmark}[1]{}
\renewcommand{\Im}{\operatorname{Im}}

\newcommand{\hol} {holomorphic}
\newcommand{\qc} {quasiconformal}

\newcommand{\Te} {Teichm\"{u}ller}

\newcommand{\iy} {\infty}
\newcommand{\fc} {\frac}

\usepackage{amsfonts}
\newcommand{\field}[1]{\mathbb{#1}}
\newcommand{\g}{\gamma}

\newcommand{\dl}{\delta}
\newcommand{\D}{\field D}
\newcommand{\om}{\omega}
\newcommand{\z}{\zeta}
\newcommand{\ov}{\overline}
\newcommand{\vp}{\varphi}
\newcommand{\hC}{\wh{\field{C}}}
\newcommand{\C}{\field{C}}
\newcommand{\R}{\field{R}}
\newcommand{\hR}{\wh{\field{R}}}
\newcommand{\B}{\mathbf{B}}
\newcommand{\T}{\mathbf{T}}

\newcommand{\Belt} {\operatorname{Belt}}

\newcommand{\uTs} {universal Teichm\"{u}ller space}

\newcommand{\const}{\operatorname{const}}

\newcommand{\Om} {\Omega}
\newcommand{\vk} {\varkappa}
\newcommand{\x} {\mathbf x}
\renewcommand{\a} {\alpha}

\newcommand{\ld}{\lambda}

\begin{document}

\title{Quasiconformal features and Fredholm eigenvalues of convex polygons}
\author{Samuel L. Krushkal}

\begin{abstract} An important open problem in geometric 
complex analysis is to find algorithms for  
explicit determination of  basic functionals intrinsically connected with conformal 
and quasiconformal maps, such as their Teichm\"{u}ller and Grunsky norms, Fredholm eigenvalues and the quasireflection coefficient. This has not been  solved even for convex polygons. 
This case has intrinsic interest in view of the connection of 
such polygons with the geometry of the universal Teichm\"{u}ller space. 

We provide a new approach,  based on affine transformations of univalent functions. 
\end{abstract}

\date{\today\hskip4mm ({\tt QcConvexPol.tex})}

\maketitle

\bigskip

{\small {\textbf {2010 Mathematics Subject Classification:} Primary:
30C55, 30C62, 30F60; Secondary: 31A35, 58B15}

\medskip

}

\bigskip

\markboth{S. L. Krushkal}{Quasiconformal features and Fredholm
eigenvalues of convex polygons}\pagestyle{headings}

\bk\bk
\centerline{\bf 1. PREAMBLE AND RESULTS}

\bigskip\noindent{\bf 1. Introductory remarks}. 

The basic functionals intrinsically connected with conformal 
and quasiconformal maps such as their Teichm\"{u}ller and Grunsky norms, the first Fredholm eigenvalue, the quasireflection coefficient imply a deep quantitative characterization of the features of these maps. 
Thus the problem to find the algorithms for  
explicit determination of these quantities is very important  
but still remains open. 

The following general result obtained in \cite{Kr2} by applying holomorphic motions 
solves this problem for unbounded convex domains giving an explicit representation of functionals by geometric characteristics of domains. 
Let $\D = \{z: \ |z| < 1\}, \ \D^* = \{z \in \hC: \ |z| > 1\}$. 

\bigskip \noindent {\bf Theorem A}. {\it For every
unbounded convex domain $D \subset \C$ with piecewise
$C^{1+\dl}$-smooth boundary $L \ (\dl > 0)$ (and all its fractional
linear images), the equalities
 \be\label{1}
 q_L = 1/\rho_L = \vk(f) = \vk(f^*) = k(f) = k(f^*)
= 1 - |\a|
\end{equation}
hold, where $f$ and $f^*$ denote the appropriately normalized
conformal maps $\D \to D$ and $\D^* \to D^* = \hC \setminus \ov{D}$, respectively, $k(f)$
and $k(f^*)$ are the minimal dilatations of their quasiconformal extensions to $\hC, \ \vk(f)$ and $\vk(f^*)$ stand for their Grunsky norms, and $\pi |\a|$ is the opening of the least interior angle between the boundary arcs $L_j \subset L$. Here $0 < \a < 1$ if the corresponding vertex is finite and $- 1 < \a < 0$ for the angle at the vertex at infinity.

The same is true for unbounded concave domains (the complements of
convex ones) which do not contain $\infty$; for those one must
replace the last term by $|\beta| - 1$, where $\pi |\beta|$ is the
opening of the largest interior angle of $D$.}

In particular, for any closed unbounded curve $L$ with the convex interior  which is $C^{1+ \dl}$ smooth at all finite points and has at infinity the asymptotes approaching the interior angle $\pi \a < 0$, we have 
$$
q_L = 1/\rho_L = 1 - |\a|.
$$

In contrast, there are bounded convex domains even with analytic
boundaries $L$ whose conformal mapping functions have different
Grunsky and Teichm\"{u}ller norms, and therefore, $\rho_L < 1/q_L$. 

\bigskip\noindent{\bf 2. Results and applications}.   

The aim of this paper is to provide the classes of bounded convex domains, especially polygons, for which these norms are equal and give explicitly the values of the associate curve functionals $k(f), \ \vk(f), \ q_L, \ \rho_L$. 

Consider the class $\Sigma^0$ of univalent functions $F(z) = z + b_0 + b_1 z^{-1} + \dots$ mapping conformally the disk $\D^*$ into $\hC \setminus \{0\}$ and having quasiconformal extensions  to $\hC$ with $F(0) = 0$. This collection naturally relates to the universal Teichm\"{u}ller space $\T$ (the space of quasisymmetric homeomorphisms $h$ of the unit circle factorized by M\"{o}bius transformations) modeled by the Schwarzian derivatives 
$$ 
S_F = (F^{\prime\prime}/F^\prime)^\prime -
(F^{\prime\prime}/F^\prime)^2/2
$$ 
of $F \in \Sigma^0$ in $\D^*$. Their inversions $f(z) = 1/F(1/z)$ form the class $S^0$ of univalent functions $f(z) = z + \sum_1^\iy a_n z^n$ of univalent functions in the unit disk $\D$ with  quasiconformal extension to $\D^*$ preserving $z = \iy$,  and $\vk(F) = \vk(f)$. 

One of the interesting questions is whether the equality of 
Teichm\"{u}ller and Grunsky norms is preserved under the affine deformations $g^c(w) = c_1 w + c_2 \overline{w} + c_3$ with $c = c_2/c_1$ (as well as of more general maps) of quasidisks. In the case of unbounded convex domains, 
this follows from Theorem A. We establish this here for bounded 
domains $D$. More precisely, we consider the maps $g^c$ which are conformal in the complementary domain $D^*  = \hC \setminus \ov{D}$ and have in $D$ a constant quasiconformal dilatation $c$. 

\bigskip\noindent{\bf Theorem 1}. {\it For any function   
$F \in \Sigma^0$ with $\vk(F) = k(F)$ mapping the disk $\D^*$ onto the complement of a bounded domain (quasidisk) $D$ and any affine deformation $g^c$ of this domain, we have the equality} 
 \be\label{2}
\vk(g^c \circ F) = k(g^c \circ F).  
\end{equation}

\bigskip 
Theorems 1 essentially increases the set of quasicircles $L \subset \hC$ for which $\rho_L = 1/q_L$ giving simultaneously the explicit values of these curve functionals. 
Even for quadrilaterals, this fact was known until now only for some special types of them (for rectangles \cite{Kr2}, [19 - 21] and for rectilinear or circular quadrilaterals having a common tangent circle \cite{We}). 

The arguments in the proof of this theorem are extended almost straightforwardly to more general case: 

\bigskip\noindent{\bf Theorem 2}. {\it Let $F \in \Sigma^0$ and $\vk(F) = k(F)$. Let $h$ be a holomorphic map $\D \to \T$ 
without critical points in $\D$ and $h(0) = S_F$. 
Denote by $\mathbf g^c$ the univalent solution of the Schwarzian  equation $S_{\mathbf g} = h(c)$ on the domain $F(\D^*)$. 
Then, for any $c \in \D$, the composition $\mathbf g^c \circ F$ also satisfies $\vk(\mathbf g^c \circ F) = k(\mathbf g^c \circ F)$.} 

\bk
Note that by the lambda lemma for holomorphic motions, the map  $h$ determines a holomorphic disk in the ball of Beltrami coefficients on $F(\D)$, which yields, together with  assumptions of the theorem, that for small $|c|$ 
$$ 
\mathbf g^c(w) = w + b_0^c + b_1^c w^{-1} + \dots \quad \text{as} \ \ w \to \iy  
$$ 
with $b_1^c \ne 0$. This is an essential point in the proof. 

The case of bounded convex polygons has an intrinsic interest, 
in view of the following negative fact underlying the features 
and contrasting Theorem A.

\bigskip\noindent{\bf Theorem 3}. {\it There exist bounded rectilinear convex polygons $P_n$ with sufficiently large number of sides such that}
 \be \label{3}
\rho_{\partial P_n} < 1/q_{\partial P_n}.
\end{equation} 

\bk
It follows simply from Theorem 1 that if a polygon $P_n$, whose edges are quasiconformal arcs, satisfies $\rho_{\partial P_n} = 1/q_{\partial P_n}$ then this equality is preserved for all its affine images. 
In particular, this is valid for all rectilinear polygons obtained by affine maps from polygons with edges having a common tangent ellipse (which includes the regular $n$-gons).      

Theorem 3 naturally gives raise to the question whether the property $\rho_{\partial P_n} = 1/q_{\partial P_n}$ is valid for all bounded convex polygons with sufficiently small number of sides. 

In the case of triangles this immediately follows from Theorem 1 as well as from Werner's result. 

Noting that the affinity preserves parallelism and moves the lines to lines, one concludes from Theorem 1 that the equality  $\rho_{\partial P_4} = 1/q_{\partial P_4}$ holds in particular for quadrilaterals $P_4$ obtained by affine transformations from quadrilaterals which are symmetric with respect to one of diagonals and for quadrilaterals whose sides have common tangent outwardly ellipse (in particular, for all parallelograms and trapezoids). 
For the same reasons, it holds also for hexagons with axial symmetry having two opposite sides parallel to this axes. 

In fact, Theorem 1 allows us to establish much stronger result 
answering the question positively for quadrilaterals.   

\bigskip\noindent{\bf Theorem 4}. {\it For every rectilinear convex quadrilateral $P_4$, we have 
 \be \label{4}
\vk(F) = k(F) = \rho_{\partial P_4} = 1/q_{\partial P_4}, 
\end{equation} 
where $F$ is the appropriately normalized conformal map of $\D^*$ onto $P_4^*$.}

\bigskip\bigskip
\centerline{\bf 2. BACKGROUND} 

\bk 
We present briefly the needed notions and results underlying the above theorems adapting those to our case; for details see, e.g. \cite{Di}, \cite{EKK}, \cite{GL}, \cite{Kr4}, \cite{Ku3}.  

\bigskip\noindent{\bf 1. A glimpse at Grunsky inequalities and
Fredholm eigenvalues}. 
Denote by $\Belt(\D)$ the unit ball of Beltrami coefficients $\mu$ supported on $\D$ and extended by zero to $\D^*$, i.e., 
$$
\Belt(\D) = \{\mu \in L_\infty(\C): \ \ \mu(z)|\D^* = 0, \ \ 
\|\mu\|_\infty < 1\}
$$
and by $w^\mu$ the solutions of the Beltrami equation 
$\partial_{\overline{z}} w = \mu \partial_z w$ on $\C$ 
with the expansion $w(z) = z + b_0 + b_1 z^{-1} + \dots$ 
in $\D^*$.  

The fundamental Grunsky theorem (extended to multiply connected
domains by Milin \cite{Mi}) states that a holomorphic function $F(z) = z + \const + O(1/z)$ in a neighborhood $U_0$ of the infinite point is extended to a univalent function on the disk $D^*$ if and only if it satisfies the inequality
$$
\Big\vert \sum\limits_{m,n=1}^\iy \ \sqrt{mn} \ \a_{m n} x_m x_n
\Big\vert \le 1,
$$
where the Grunsky coefficients $\a_{mn}(f)$  are determined by
$$
\log \fc{F(z) - F(\z)}{z - \z} = - \sum\limits_{m,n=1}^\iy \  \a_{mn} z^{-m} \z^{-n}, \quad (z, \z) \in (\D^*)^2,
$$
taking the principal branch of the logarithmic function, and
$\mathbf x = (x_n)$ ranges over the unit sphere $S(l^2)$ of the
Hilbert space $l^2$ of sequences with $\|\mathbf x\|^2 =
\sum\limits_1^\iy |x_n|^2$ (cf. \cite{Gr}). The quantity
$$
\vk(F) = \sup \Big\{\Big\vert \sum\limits_{m,n=1}^\iy \ \sqrt{mn} \
\a_{mn} x_m x_n \Big\vert: \mathbf x = (x_n) \in S(l^2) \Big\} 
$$
is called the {\it Grunsky norm} of the map $F$.

It is dominated by the {\it Teichm\"{u}ller norm} $k(F)$ of 
this map, i.e., with the minimal dilatation among quasiconformal extensions of $F$ onto $\D$ (see \cite{Ku1}, \cite{Kr8}); so,
  \be\label{5}
\vk(F) \le k(F) = \tanh \tau_\T(\mathbf 0, S_F),
\end{equation}
where $\tau_\T$ denotes the Teichm\"{u}ller distance on $\T$).  The second norm is intrinsically connected with integrable
holomorphic quadratic differentials on $\D$ (the elements of the subspace $A_1 = A_1(\D)$ of $L_1(\D)$ formed by holomorphic
functions), while the Grunsky norm naturally relates to the {\it abelian} structure determined by the set of quadratic differentials
$$
A_1^2 = \{\psi \in A_1: \ \psi = \om^2\};
$$
having only zeros of even order on $\D$. In terms of the pairing
$$
\langle \mu, \psi\rangle_\D = \iint\limits_\D  \mu(z) \psi(z) dx dy,
\quad \mu \in L_\iy(\D), \ \psi \in L_1(\D) \ \ (z = x + iy),
$$
we have the following results characterizing the functions with
$\vk(F) = k(F)$.

\bigskip\noindent
{\bf Lemma 1}. \cite{Kr1}, \cite{Kr8} \ {\it For all $F = F^\mu \in \Sigma^0$,
$$
\vk(F) \le k \fc{k + \a(F)}{1 + \a(F) k}, \quad k = k(F),
$$
and $\vk(F) < k$ unless
 \be\label{6}
\a(F) := \sup \ \{|\langle \mu, \psi\rangle_\D|: \ \psi \in A_1^2, \ \|\psi\|_{A_1(\D)} = 1\} = \|\mu\|_\iy;
\end{equation}
the last equality is equivalent to $\vk(F) = k(F)$. Moreover, for small $\|\mu\|_\iy$,
$$
\vk(F) = \sup \ |\langle \mu, \psi\rangle_\D| +
O(\|\mu\|_\iy^2), \quad \|\mu\|_\iy \to 0,
$$
with the same supremum as in (6).

If $\vk(F) = k(F)$ and the equivalence class of $F$ (the collection of maps equal to $F$ on $S^1 = \partial D^*$) is a Strebel point, then the extremal $\mu_0$ in this class is necessarily of the form }
 \be\label{7}
\mu_0 = \|\mu_0\|_\iy |\psi_0|/\psi_0 \ \ \text{with} \ \ \psi_0 \in
A_1^2.
\end{equation}

Geometrically, (6) means the equality of the Carath\'{e}odory and Teichm\"{u}ller distances on the image of the geodesic disk
$\D(\mu_0) = \{t\mu_0 /\|\mu_0\|_\iy: \ t \in \D\}$ in the space
$\T$. For functions $F \in \Sigma^0$ holomorphic in the closed disk $\ov{\D^*}$, the relation (7) was also obtained by a different method in \cite{Ku4}. 

An important property of the Grunsky coefficients $\a_{mn}(F) = \a_{mn}(S_F)$ is that these coefficients are holomorphic 
functions of the Schwarzians $\vp = S_F$ on the \uTs \ $\T$.
Therefore, for every $F \in \Sigma^0$ and each $\x = (x_n) \in S(l^2)$, the series
 \be\label{8}
h_{\mathbf x}(\vp) =
\sum\limits_{m,n=1}^\iy \ \sqrt{mn} \ \a_{mn}(\vp) x_m x_n
\end{equation}
defines a holomorphic map of the space $\T$ into the unit disk $\D$, and 
$\vk(F) = \sup_{\x} | h_{\mathbf x}(S_F)|$. 

The convergence and holomorphy of the series (8) simply
follow from the inequalities
$$
\Big\vert \sum\limits_{m=j}^M \sum\limits_{n=l}^N \ \sqrt{m n} \
\a_{mn} x_m x_n \Big\vert^2 \le \sum\limits_{m=j}^M |x_m|^2
\sum\limits_{n=l}^N |x_n|^2
$$
(for any finite $M, \ N$) which, in turn, are a consequence of the classical area theorem (see, e.g., [24, p. 61]).
  
Using Parseval's equality, one obtains that the elements of the distinguished set $A_1^2$ are represented in the form
 \be\label{9}
\psi(z) = \fc{1}{\pi} \sum\limits_{m+n=2}^{\iy}
\sqrt{mn} \ x_m x_n z^{m+n-2}
\end{equation}
with $\mathbf x = (x_n) \in l^2$ so that $\|\mathbf x\|_{l^2} = \|\psi\|_{A_1}$ (see \cite{Kr1}). 

\bigskip
A crucial point here is that for a generic function $F \in
\Sigma^0$ in (5) the strict inequality $\vk(F) < k(F)$ is valid; moreover, it holds on the (open) dense subset of $\Sigma^0$ in both strong and weak topologies (i.e., in the Teichm\"{u}ller
distance and in locally uniform convergence on $D^*$); see
\cite{Kr1}, \cite{Kr6}, \cite{KK}, \cite{Ku2}, \cite{Ku3}. So it is important to know whether for a concrete function $F$, we have $\vk(F) = k(F)$ . This fact is deeply related 
to various topics in the complex geometry of the Teichm\"{u}ller space theory, geometric complex analysis, Fredholm eigenvalues and boundary problems, operator theory, etc.

\bigskip\noindent{\bf 2. Quasireflections}. The 
quasiconformal reflections (or quasireflections)
represent a special case of topological orientation reversing
involutions of the sphere $S^2 = \hC = \C \cup \C$. Any
quasireflection preserves pointwise fixed a quasicircle $L \subset \hC$ interchanging its inner and outer domains (because, due to \cite{Kr4}, any set $E \subset S^2$, which admits quasireflections, is necessarily located on a quasicircle with the same reflection coefficient).

One defines for each mirror $E$ its {\it reflection coefficient}
$q_E = \inf \| \partial_z f/\partial_{\ov z} f \|_\iy$ (taking the
infimum over all quasireflecions across $E$) and quasiconformal
dilatation  $Q_E = (1 + q_E)/(1 - q_E) \geq 1$. Due to \cite{Ah2}, \cite{Kr4}, \cite{Ku4}, 
$$
Q_E = (1 + k_E)^2/(1 - k_E)^2,
$$
where $k_E = \inf \|\partial_{\ov z} f/\partial_z f \|_\iy$ over all quasicircles $L \supset E$ and all orientation preserving
quasiconformal homeomorphisms $f: \ \hC \to \hC$ with $f(\hR) = L$.

\bigskip\noindent{\bf 3. Fredholm eigenvalues}. 
The {\it Fredholm eigenvalues} $\rho_n$ of a smooth closed Jordan
curve $L \subset\hC$ are the eigenvalues of its double-layer
potential, or equivalently, of the integral equation
$$
u(z) +  \fc{\rho}{\pi} \int\limits_L \ u(\z) \fc{\partial}{\partial
n_\z} \log \frac{1}{|\z - z|} ds_\z = h(z),
$$
which has has many applications (here $n_\z$ is the outer normal and $ds_\z$ is the length element at $\z \in L$). 

The least positive eigenvalue $\rho_L = \rho_1$ plays a crucial role
and is naturally connected with conformal and quasiconformal maps .
It can be defined for any oriented closed Jordan curve $L$ by
$$
\fc{1}{\rho_L} = \sup \ \fc{|\mathcal D_G (u) - \mathcal D_{G^*}
(u)|} {\mathcal D_G (u) + \mathcal D_{G^*} (u)},
$$
where $G$ and $G^*$ are, respectively, the interior and exterior of
$L; \ \mathcal D$ denotes the Dirichlet integral, and the supremum
is taken over all functions $u$ continuous on $\hC$ and harmonic on
$G \cup G^*$. In particular, $\rho_L = \iy$ only for the circle.

An upper bound for $\rho_L$ is given by Ahlfors' inequality
 \be\label{10}
\fc{1}{\rho_L} \le q_L,
\end{equation}  
where $q_L$ denotes the minimal dilatation of quasireflections
across $L$.  This inequality is equivalent to (4) and serves as a background for defining the value $\rho_L$. This value is intrinsically connected with the Grunsky
operator, which is qualitatively expressed by the 
K\"{u}hnau-Schiffer theorem; it states that
 $\rho_L$ is reciprocal to the Grunsky norm $\vk(f)$ of the Riemann mapping function of the exterior domain of $L$ (cf. \cite{Ku2}, \cite{Sc}).

\bigskip\noindent{\bf 4. Plurisubharmonicity of the Teichm\"{u}ller metric}. 
Due to the fundamental Gardiner-Royden theorem, the Kobayashi  and Teichm\"{u}ller metrics on Teichm\"{u}ller spaces are equal.  
An essential strengthening of this theorem for the space $\T$
established in \cite{Kr3} by applying the Grunsky coefficient technique  yields 

\bk\noindent{\bf Lemma 2}. \cite{Kr3} {\it The infinitesimal forms of both metrics $\mathcal K_\T(\vp, v)$ and $F_\T(\vp, v)$ on the tangent bundle $\mathcal T(\T)$ of $\T$ are continuous logarithmically plurisubharmonic in $\vp \in \T$ and have constant holomorphic sectional curvature $\kappa_{\mathcal K}(\vp, v) = - 4$.} 

In addition, these infinitesimal metrics are Lipshitz continuous (see \cite{EE}). 
The global distances (integrated forms of these metrcis) are logarithmically plurisubharmonic in each of their variables on $\T \times \T$ (cf. \cite{Kr4}).

\bigskip\bigskip
\centerline{\bf 3. PROOF OF THEOREM 1}. 

\bigskip 
The proofs of all the above theorems essentially rely 
on the following deep fact stated as a conjecture in \cite{KK} and proven in \cite{Kr5}.

\bigskip\noindent{\bf Proposition 1}. {\it Any sequence of the
functions $F_n \in \Sigma^0$ with $\vk(F_n) = k(F_n)$ cannot
converge locally uniformly in $\D^*$ to a function $F \in \Sigma^0$ with $\vk(F) < k(F)$. }

Thus it suffices to establish the assertion of the theorem for functions $F \in \Sigma^0$ which are holomorphic on the closed disk $\ov{\D^*}$. Indeed, by the density theorem of \cite{Kr5}, the Strebel points $\vp = S_f^{k|\mu|/\mu} \in \T$ with $\mu \in A_1^2$ representing the functions $F \in \Sigma^0$ with equal Grunsky and Teichm\"{u}ller norms are dense in $\B$ norm in the set of all $F \in \Sigma^0$ with $\vk(F) = k(F)$. So one can pass to $F^{\mu_r}$ with $\mu_r(z) = k |\mu(rz)|/\mu(rz)$ taking $r < 1$ so that 
the homotopy disk $\D(S_{F^{\mu_r}})$ has no critical points 
in the annulus $\{r < |t| < 1\}$. Then the equality (2) in the limit $r \to 1$ follows again by applying Proposition 1.  
We accomplish the proof of the theorem in three stages. 

\bk\noindent
$\mathbf{1^0}$. \ First, we establish some auxiliary results characterizing the homotopy disk of a map with $\vk(F) = k(F)$. 

Take the generic homotopy function  
$$
F_t (z) = t F(z/t) =  z + b_0 t + b_1 t^2 z^{-1} + b_2 t^3 z^{-2} + \dots: \ \D^* \times \D \to \hC.   
$$ 
Then $S_{F_t}(z) = t^{-2} S_F(t^{-1} z)$ and this point-wise map determines a \hol \ map
$\chi_F(t) =  S_{F_t}(\cdot): \ \D \to \T$ so that the homotopy disks $\D(S_F) = \chi_F(\D)$ have only cuspidal critical points and foliate the space $\T$. Note also that 
$$
\a_{m n}(F_t) = \a_{m n}(F) t^{m+n},  
$$ 
and if $f(z) = 1/F(1/z)$ maps the unit disk onto a convex domain, then all level lines $F(|z| = r)$ for $z \in \D^*$ are starlike.

\bk\noindent{\bf Lemma 3}. {\it If the homotopy function $F_t$ of $F \in \Sigma^0$ satisfy $\vk(F_{t_0}) = k(F_{t_0})$ for some 
$0 < t_0 < 1$, then the equality $\vk(F_t) = k(F_t)$ holds 
for all $|t| \le t_0$ and the homotopy disk $\D(S_{f_t})$ has 
no critical points $t$ with $0 < |t| < t_0$.} 

\bigskip\noindent{\bf Proof}. Take the univalent extension 
$F_1$ of $F$ to a maximal disk $\D_b^* = \{z \in \hC: \ |z| > b\}, \ (0 < b < 1)$ and define  
$$ 
F^*(z) = b^{-1} F_1(b z) \in \Sigma^0, \quad |z| > 1.
$$ 
Its Beltrami coefficient in $\D$ is defined by holomorphic quadratic differentials $\psi \in A_1^2$ of the form (9),  
and we have the holomorphic map 
 \be\label{11}
h_{\mathbf x^b}(S_{F_t^*}) =
\sum\limits_{m,n=1}^\iy \ \sqrt{mn} \ \a_{mn}(F^*) x_m^b x_n^b 
(b t)^{m+n}   
\end{equation}  
of the disk $\D(S_{F^*})$ into $\D$. In view of our assumption on $F$, the series (11) is convergent in some wider disk $\{|t| < a\} (a > 1)$. 

Using the map (11), we pull back the hyperbolic metric 
$\ld_\D(t) = |dt|/(1 - |t|^2)$ to the disk $\D(S_{F_1})$ (parametrized by $t$) and define on this disk 
the conformal metric $ds = \ld_{\wt h_\x}(t) |dt|$ with 
 \be\label{12}
\ld_{\wt h_{\x^b}}(t) = (h_{\x^a} \circ \chi_{f_1})^* \ld_\D = 
\fc{|\wt h_{\x^b}\prime (t)| |dt|}{1 - |\wt h_{\x^b}(t)|^2}.
\end{equation} 
of Gaussian curvature $- 4$ at noncritical points. In fact, this is the supporting metric at $t = a$ for the upper envelope 
$\ld_\vk = \sup_{\x \in S(l^2)} \ld_{\wt h_{\x^b}}(t)$ of metrics (12) followed by its upper semicontinuous regularization 
$u(t) =\limsup_{t^\prime\to t} u(t^\prime)$ 
(supporting means that $\ld_{\wt h_{\x^b}}(a) = \ld_\vk(a)$  
and $\ld_{\wt h_{\x^b}}(t) < \ld_\vk(t)$ in a neighborhood of $a$). 
  
The metric $\ld_\vk(t)$ is logarithmically subharmonic on $\D$  and its generalized Laplacian  
$$
\Delta u(t) = 4 \liminf\limits_{r \to 0} \frac{1}{r^2}
\Big\{ \frac{1}{2 \pi} \int_0^{2\pi} u(t + re^{i \theta}) d
\theta - \ld(t) \Big\} 
$$ 
satisfies 
$$
\Delta \log \ld_\vk \ge 4 \ld_\vk^2   
$$ 
(while for $\ld_{\wt h_{\x^b}}$ we have at its noncritical points $\Delta \log \ld_{\wt h_{\x^b}} =  4 \ld_{\wt h_{\x^b}}^2$). 

Note also that the Grunsky coefficients define on the tangent bundle $\mathcal T(\T)$ a new Finsler structure $F_\vk(\vp, v)$  dominated by the infinitesimal Teichm\"{u}ller metric $F(\vp, v)$. This structure generates on any embedded \hol \ disk $\g(\D) \subset \T$ the corresponding Finsler metric $\ld_\g(t) = F_\vk(\g(t), \g'(t))$ and reconstructs the Grunsky norm by integration along the Teichm\"{u}ller disks: 

\bigskip\noindent
{\bf Lemma 4}. \cite{Kr3} {\it On any extremal Teichm\"{u}ller disk $\D(\mu_0) = \{\phi_\T(t \mu_0): \ t \in\D\}$
(and its isometric images in $\T$), we have the equality}
$$
\tanh^{-1}[\vk(f^{r\mu_0})] = \int\limits_0^r
\ld_\vk(t) dt.
$$

Taking into account that the disk $\D(S_f)$ touches at the point
$\vp = S_{f_a}$  the Teichm\"{u}ller disk centered at the origin of $\T$ and passing through this point and that the metric $\ld_\vk$ does not depend on the tangent unit vectors whose initial points are the points of $\D(S_f)$, one obtains from Lemma 3 and the equality 
$\vk(f_a) = k(f_a)$ that also 
 \be\label{13}
\ld_\vk(a) = \ld_{\mathcal K}(a). 
\end{equation} 

We compare the metric $\ld_{\wt h_{\x^b}}$ with $\ld_{\mathcal K}$ using Lemma 2 and Minda's maximum principle given by 

\bigskip\noindent{\bf Lemma 5}. \cite{Min} {\it If a function
$u : \ D \to [- \iy, + \iy)$ is upper semicontinuous in a domain $D \subset \C$ and its (generalized) Laplacian satisfies the inequality $\Delta u(z) \ge K u(z)$ with some positive constant $K$ at any point $z \in D$, where $u(z) > - \iy$, and if
$$
\limsup\limits_{z \to \z} u(z) \le 0 \ \ \text{for all} \
\z \in \partial D,
$$
then either $u(z) < 0$ for all $z \in D$ or else $u(z) = 0$ for all $z \in D$.} 

Take a sufficiently small neighborhood $U_0$ of the point $t = a$, and let 
$$
M = \{\sup \ld_{\mathcal K}(t): t \in U_0\};
$$
then in this neighborhood,
$\ld_{\mathcal K}(t) + \ld_{\wt h_{\x^b}}(t) \le 2M$ 
and the function
$$
u = \log \fc{\ld_{\wt h_{\x^b}}}{\ld_{\mathcal K}} = 
\log \ld_{\wt h_{\x^b}} - \log \ld_{\mathcal K}   
$$
satisfies 
$$
\Delta u = 4 (\ld_{\wt h_{\x^b}}^2 - \ld_{\mathcal K}^2) \ge
8M (\ld_{\wt h_{\x^b}} - \ld_{\mathcal K}).
$$
The elementary estimate
$$
M \log(t/s) \ge t - s \quad \text{for} \ \ 0 < s \le t < M
$$
(with equality only for $t = s$) implies that
$$
M \log \fc{\ld_{\wt h_{\x^b}}(t)}{\ld_{\mathcal K}(t)} \ge \ld_{\wt h_{\x^b}}(t) - \ld_{\mathcal K}(t),
$$
and hence,
$$
\mathbf{\D} u(t) \ge 4 M^2 u(t).
$$
Lemma 5 and the equality (13) imply that the metrics 
$\ld_{\wt h_{\x^b}}, \$\ld_\vk, \ld_{\mathcal K}$ must be equal in the entire disk $\D(S_F)$, which yields by Lemma 3 the equality 
$$
\vk(F_r) = k(F_r) = \Big\vert\sum\limits_{m,n=1}^\iy \ \sqrt{mn} \ \a_{mn}(F_1)
r^{m+n} x_m^r x_n^r \Big\vert 
$$
for all $r = |t| \in (0, 1)$ (with $(x_n^r) \in S(l^2)$ depending on $r$) and that for any $F \in \Sigma^0$ 
with $\vk(F) = k(F)$ {\it its homotopy disk $\D(S_F)$ has only a singularity at the origin of $\T$}. 

\bk\noindent
$\mathbf{2^0}$. \ We may now investigate the action of affine deformations on the set of functions $F \in \Sigma^0$ with equal Grunsky and Teichm\"{u}ller norms.    

\noindent{\bf Lemma 6}. {\it For any affine deformation $g^c$ of a convex domain $D$ with expansion 
$g^c(w) = w + b_0^c + b_1^c w^{-1} + \dots$ near $w = \iy$, 
we have  
 \begin{equation}\label{14}
b_1^c = \fc{S_{g^c}(\iy)}{6} = \fc{1}{6} \lim\limits_{z\to \iy} 
w^4 S_{g^c}(w) \ne 0,  
\end{equation} 
and for sufficiently small $|c|$ all composite maps   
 \begin{equation}\label{15}
W_{F,c}(z) = g^c \circ F(z) = z + \wh b_0^c + \wh b_1^c z^{-1} + \dots, 
\quad F \in \Sigma^0
\end{equation}  
also satisfy $\wh b_1^c \ne 0$. }

\bk\noindent{\bf Proof}. One can assume that $\G = \partial D$ is a smooth curve. For small $|c|$, 
$$
g^c(w) = w - \fc{c}{\pi} \iint\limits_D \fc{dx d y}{z - w} + O(c^2),  
$$
and by the Cauchy-Green formula,  
$$
\iint\limits_D \fc{dx d y}{z - w} = 2 i \int\limits_\G \fc{\ov{z} dz}{z - w} = - \fc{2i}{w} \int\limits_\G \ov{z} dz 
+ O\Bigl(\frac{1}{w^2}\Bigr) = - \frac{1}{w}  \iint\limits_D dx dy +  O\Bigl(\fc{1}{w^2}\Bigr), \quad w \to \infty; 
$$ 
hence, $b_1 \ne 0$, which proves (14). 

The second assertion of the lemma follows from (14) and the equality 
$$
6 \wh b_1^c = \lim\limits_{z\to \iy} z^4 S_{W_{F,c}}(z) = 
\lim\limits_{z\to \iy} z^4 \left[(S_{g^c}\circ F)(z) (F^\prime(z)^2 + S_F(z)\right] 
$$
for any $F \in \Sigma^0$, which completes the proof. 

Lemma 6 allows one to apply to compositions $W_{F,c}$  the following result of K\"{u}hnau \cite{Ku2}: 

\bigskip\noindent{\bf Lemma 7}. {\it For any function $F(z) = z + b_0 + b_1 z^{-1} + \dots \in \Sigma(0)$ with $b_1 \ne 0$, the extremal quasiconformal extensions of the homotopy functions $F_t$ to $\D$ are defined for sufficiently small $|t| \le r_0 =  r_0(F) \ (r_0 > 0)$ by nonvanishing holomorphic quadratic differentials, and therefore, $\vk(F_t) = k(F_t)$}. 

\bk
It follows from Lemmas 6 and 7 that for any $F \in \Sigma^0$ and any affine transformation $g^c$ of domain $F(\D)$ the homotopy functions 
$$
W_{F,c,t} = g^c \circ F_t
$$
of maps (15) satisfy 
$$
\vk(g^c \circ F_t) = k(g^c \circ F_t) \quad \text{for all} \ \ 
|t| \le r_0(F, g^c) \ \ (r_0(F, g^c)).    
$$ 

Note also that since $F$ was chosen to be holomorphic on $\ov{\D^*}$ and its homotopy disk has only a singularity at $t = 0$, the equality 
$$
S_{W_{F,c}}(z) = (S_{g^c}\circ F)(z) (F^\prime(z)^2 + S_F(z) 
$$
implies that for sufficiently small $|c|$ the homotopy disk 
$\D(S_{W_{F_t,c}})$ of $W_{F_t,c}$ also has only a singularity at $t = 0$. 
Fix a such $c$.  

Using the restrictions to $\D(S_{W_{F_t,c}})$  
of the corresponding holomorphic maps 
 \begin{equation}\label{16}
h_\x(S_{W_{F,c}}) = \sum\limits_{m,n=1}^\iy \ \sqrt{mn} \ 
\a_{m n}(W_{F,c}) x_m x_n: \ \T \to \D
\end{equation} 
one obtains for 
$\wh h_x = h_\x(S_{W_{F,c}}) \circ \chi_{W_{F,c}}$, 
the conformal metrics 
$$
\ld_{\wh h_\x}(t) = \wh h_x^* \ld_\D = 
\fc{|\wt h_\x^\prime(t)| |dt|}{1 - |\wt h_\x(t)|^2} 
$$ 
of curvature $- 4$. 
Their upper envelope  $\ld_\vk(t)$  
(followed by its upper semicontinuous regularization) is a  subharmonic metric  of generalized Gaussian curvature $\kappa_{\ld_\vk} \le -4$. 
It can be compared with the infinitesimal Kobayashi metric 
$\ld_{\mathcal K}$ on $\D(S_{W_{F_t,c}})$ similar to the previous step $\mathbf{1^0}$ by applying Lemmas 3 and 4. 
This implies the equalities $\ld_\vk = \ld_{\mathcal K}$ on $\D(S_{W_{F,c}})$ and 
 \begin{equation}\label{17}
\vk(W_{F,c}) = k(W_{F,c}). 
\end{equation} 

\bk\noindent 
$\mathbf{3^0}$. \ It remains to extend the last equality to all 
$c$ with $|c| < 1$. 
Noting that by the chain rule for Beltrami coefficients $\mu, \nu$ from the unit ball in $L_\iy(\C)$,  
$$
w^\mu \circ w^\nu = w^\tau \ \ \text{with} \ \ \tau = (\nu + \wt{\mu})/(1 + \ov{\nu} \wt{\mu}) 
$$
and $\wt{\mu}(z) = \mu(w^\nu(z)) \ov{w^\nu_z}/w^\nu_z$ (so    
for $\nu$ fixed, $\tau$ depends holomorphically on $\mu$ in $L_\iy$ norm), one can regard (16) as holomorphic functions of $c \in \D$ and construct in a similar way the corresponding Finsler metrics 
$$
\ld_{h_\x}(c) = |\wt h_\x^\prime (c)| |dc|/(1 - |\wt h_\x(c)|^2), \quad |c| < 1. 
$$
Now take the upper envelope $\ld_\vk(c)$ of these metrics and its upper semicontinuous regularization getting now a subharmonic metric of Gaussian curvature $\kappa_{\ld_\vk} \le -4$ on the nonsingular disk $\{|c| < 1\}$. One can repeat for 
this metric all the above arguments using the already established equality (17) for small $|c|$. The assertion of Theorem 1 follows again by applying Lemmas 4 and 5, which completes the proof of Theorem 1.

\bigskip\bigskip 
\centerline{\bf 4. PROOF OF THEOREM 3} 

\bk
Take a function $F_0(z) = z + b_0 + b_m z^{-m} + \dots \in \Sigma^0$ with $m > 1, \ F_0(0) = 0$, mapping the disk $\D^*$ onto a domain $D$ with smooth boundary, hence having Teichm\"{u}ller quasiconformal extension to $\hC$ with Beltrami coefficient
$$
\mu_0(z) = k_0 |\psi(z)|/\psi(z), \quad |z| < 1,
$$
where $\psi_0 \in A_1(\D)$ is of the form
$$
\psi_0(z) = c_p z^p + c_{p+1}z^{p+1} + \dots \quad  c_p \ne 0, \quad p \ \text{odd} 
$$
(then $\vk(F_0) < k(F_0) = k_0$) and such that also its homotopy functions
$$
F_{0t}(z) = t F_0(z/t), \quad z \in \D^*, \ \ 0 \le t < 1,
$$
satisfy
 \begin{equation}\label{18}
\vk(F_{0t}) < k(F_{0t}) \quad \text{for} \ \ 0 < t < 1
\end{equation}
(one can use, for example, the map $F_0(z)$ conformal in $\D^*$ and having on the unit disk the Beltrami coefficient $\mu_0(z) = k_0
|z|^p/z^p$ with odd $p$). Its inversion $f_0(z) = 1/F_0(1/z) = z - b_0 z^2 + \dots$ maps, by the well-known geometric
properties of univalent functions, the disks $\D_r = \{|z| < r\}, \ r < 1$,   images $F_0(\D_r)$ of the disks
$\D_r = \{|z| < r\}$ onto convex domains for all $r \le 2 -
\sqrt{3}$, and 
$$
\vk(f_0) = \vk(F_0) < k(F_0) = k(f_0). 
$$ 
Take a fixed $r < 2 - \sqrt{3}$ so that $\|S_{f_{0r}}\|_\B < 2$ and a dense subset $E = \{z_1, z_2, \dots, z_n, \dots\}$ on the unit circle $S^1 = \partial \D$. Now consider the convex rectilinear polygons $P_n$ located in the interior of quasicircle $L_0 = f_{0r}(S^1)$ with vertices
$f_{0r}(z_1), \dots, f_{0r}(z_n)$ on $L_0$, and let $F_{P_n}$ be an appropriately normalized conformal map of $\D^*$ onto the
complement of $P_n$. Then, by (18) and Proposition 1, there exists a natural $n_0$ such that 
$$
\vk(F_{P_n}) < k(F_{P_n})
$$ 
for any $P_n$ with $n > n_0$. In view of invariance of both sides of (10) under the M\"{o}bius maps of $\hC$ and by the K\"{u}hnau-Schiffer theorem the last inequality is equivalent to (3). This completes the proof of Theorem 3.

\bk
\bigskip\noindent
\centerline{\bf 5. PROOF OF THEOREM 4} 
 
\bk 
In view of Proposition 1, it sufficed to prove the theorem 
for bounded quadrilaterals $P_4 = A_1 A_2 A_3 A_4$ with vertices $A_j$ (ordered according to positive direction of $\partial P_4$) such that that the line in $P_4$ drawn from the vertex $A_1$ parallel to the opposite edge $A_2 A_3$ separates 
this edge from the remaining vertex $A_4$. 

Fix such a quadrilateral $P_4^0 = A_1^0 A_2^0 A_3^0 A_4^0$ and 
consider the collection $\mathcal P^0$ of quadrilaterals $P_4 = A_1^0 A_2^0 A_3^0 A_4$ with the same first three vertices and variable $A_4$; the corresponding $A_4$ runs over a subset $E$  of the trice punctured sphere $\hC \setminus \{A_1^0, A_2^0, A_3^0\}$.

The conformal map $F$ of the disk $\D^*$ onto the complementary domain $P_4^* = \hC \setminus \ov{P_4}$ is represented by the Schwarz-Christoffel integral
 \begin{equation}\label{19}
F(z) = d_1 \int\limits_0^z \prod\limits_1^4(\z - e_j)^{\a_j - 1} \ \frac{d \z}{\z^2} + d_0,
\end{equation}
where $e_j = F^{-1}(A_j) \in S^1, \ \pi \a_j$ is the interior 
angle at $A_j$ for $P_4^*$, and $d_0, d_1$ are two complex constants. Let $F^0$ denote the conformal map for the complement of $P_4^0$. 

One obtains from the general properties of quasiconformal maps and (19) that the logarithmic derivatives 
$b_F = (\log F^\prime)^\prime = F^{\prime\prime}/F^\prime$
of maps $F$ defining the quadrilaterals $P_4 \in \mathcal P^0$ are (for a fixed $z$)  real analytic functions of $t = A_4$. Passing to their Schwarzians 
$$ 
S_F = b_F^\prime - \fc{1}{2} b_F^2 \in \T    
$$
one can find a smooth real arc $\G = \mathbf b(t) \subset \T$ containing the point $S_{F^0}$ and the points corresponding to trapezoids; here $\mathbf b$ denotes the map $t = A_4 \to S_F$. 

Since $\T$ is a domain, there is a tubular neighborhood 
containing $\G$ therefore, $\G$ is located on some nonsingular \hol \ disk of the form 
$\Om_0 = \mathbf F(G_0) \subset \T$, where $G_0$ is a simply connected planar domain containing the set $E$. This disk is not geodesic in the Teichm\"{u}ller-Kobayashi metric on $\T$ and does not pass through the basepoint $\vp = \mathbf 0$ of this space, but one can apply to it the same arguments as in the proof of Theorem 1 constructing similar to (8) the holomorphic maps
$$ 
h_{\x}(\vp) = \sum\limits_{m,n=1}^\iy \ \sqrt{m n} \
\a_{m n}(F) x_m x_n : \ \T \to \D \ \ (\vp = S_F, \ \ \x \in S(l^2)). 
$$
The restrictions of these maps to the disk $\Om_0$ determine (again by pulling back the hyperbolic metric of the disk) the corresponding conformal metrics 
$$
\ld_{\wh h_\x}(t) =  
\fc{|\wt h_\x^\prime(t)| |dt|}{1 - |\wt h_\x(t)|^2},    
$$ 
and their upper semicontinuous envelope $\ld_\vk$ is a subharmonic metric on $\Om_0$ of (generalized) Gaussian curvature $\kappa(\ld_{\vk}) \leq - 4$.

Noting that the collection $\mathcal P^0$ contains the trapezoids, for which we have the equalities (4) by Theorem 1  
(and consequently, the infinitesimal equality (13) at the corresponding points $t$),  
one again obtains by applying Lemma 4 that the constructed metric $\ld_\vk$ must coincide at all points of $\Om_0$
with the dominant infinitesimal Teichm\"{u}ller-Kobayashi metric $\ld_{\mathcal K}$ of $\T$. Together with Lemma 5, 
this provides the global equalities (4) for all points of 
the disk $\Om_0$, which yields the assertion of Theorem 4 
for a given quadrilateral $P_4^0$.

\bk
\bigskip\noindent
\centerline{\bf 6. ADDITIONAL REMARKS}  

\bigskip\noindent
{\bf 1.} Generically in Lemma 7, $r_0(F) < 1$; this is caused by the critical points of the disk $\D(S_F)$ and circular symmetry of both infinitesimal metrics $\ld_\vk$ and $\ld_\mathcal{K}$ on this disk.   

\bigskip\noindent
{\bf 2.} Another reason why the convex polygons are interesting for quasiconformal theory is their close geometric connection with the geometry of universal Teichm\"{u}ller space.
The relations of type (1) are valid for bounded convex rectilinear polygons $P_n$ in the following truncated form. Denoting the vertices of $P_n$ by $A_j$ and their interior 
angles by $\pi \a_j \ (j = 1, \dots, n)$, one represents the conformal map $f_n$ of the upper
half-plane $H = \{z: \ \Im z > 0\}$ onto $P_n$ by the
Schwarz-Christoffel integral
$$
f_n(z) = d_1 \int\limits_0^z (\xi - a_1)^{\alpha_1 - 1} (\xi -
a_2)^{\alpha_2 - 1} ... (\xi - a_n)^{\alpha_n - 1} d \xi + d_0,
$$
(with $a_j = f_n^{-1}(A_j) \in \R$ and complex constants $d_0,
d_1$). Its Schwarzian derivative is given by
$$
S_{f_n}(z) = b_{f_n}^\prime(z) - \fc{1}{2}  b_{f_n}^2(z)
= \sum\limits_1^n \frac{C_j}{(z - a_j)^2} - \sum\limits_{j,l=1}^n
\frac{C_{jl}}{(z - a_j)(z - a_l)},
$$
where $C_j = \alpha_j - 1 - (\alpha_j - 1)^2/2 < 0, \ \ C_{jl} =
(\alpha_j - 1)(\alpha_l - 1) > 0$. It is a point of the universal
Teichm\"{u}ller space $\T$ modeled as a bounded domain in the space
$\B(H)$ of hyperbolically bounded holomorphic functions on $H$ with
norm $\|\vp\|_\B = \sup_H |z - \ov z|^2 |\vp(z)|$. 

Denote by $r_0$ the positive root of the equation
$$
\frac{1}{2} \Bigl[\sum\limits_1^n (\alpha_j - 1)^2 +
\sum\limits_{j,l=1}^n (\alpha_j - 1)(\alpha_l - 1) \Bigr] r^2 -
\sum\limits_1^n (\alpha_j - 1) \ r - 2 = 0,
$$
and put $ S_{f_n,t} = t b_{f_n}^\prime -  b_{f_n}^2/2, \ t > 0$. Then we have

\bigskip\noindent
{\bf Proposition 2}. \cite{Kr7} {\it For any convex polygon $P_n$, the Schwarzians $r S_{f_n,r_0}$  define for any $0 <r < r_0$ a univalent function $w_r: H \to \C$ whose harmonic Beltrami coefficients $\nu_r(z) = - (r/2) y^2 S_{f_n,r_0}(\overline z)$ is extremal in its equivalence class, and}
$$
k(w_r) = \vk(w_r) = \fc{r}{2} \| S_{f_n,r_0}\|_\B.
$$

By the Ahlfors-Weill theorem \cite{AW}, every $\vp \in \B(H)$ with $\|\vp\|_\B < 1/2$ is the Schwarzian derivative
$S_W$ of a univalent function $W$ in $H$, and this function has quasiconformal extension onto the lower half-plane $H^* = \{z: \ \Im z < 0\}$ with Betrami coefficient of the form
$$
\mu_\vp(z) = - 2 y^2 \vp(\ov z), \quad \vp = S_f \ (z = x + i y \in H)
$$
called harmonic.
Proposition 2 yields that any $w_r$ with $r < r_0$ does not admit extremal quasiconformal extensions of Teichm\"{u}ller type, and in view of extremality of harmonic coefficients $\mu_{S_{w_r}}$ the Schwarzians $S_{w_r}$ for some $r$ between $r_0$ and $1$ must lie outside of the space $\T$; so this space is not a starlike domain in $\B$.

\medskip
{\small\em{ \leftline{Department of Mathematics, Bar-Ilan
University} \leftline{5290002 Ramat-Gan, Israel} \leftline{and
Department of Mathematics, University of Virginia,}
\leftline{Charlottesville, VA 22904-4137, USA}}}

\end{document}